# Chaotic distribution of prime numbers and digits of π

Marek Berezowski

Cracow University of Technology, Poland

In the paper the distribution of prime numbers and digits of π were presented as chaotic.

The first analysis was based on sequence $N=100000$ of prime numbers and a graph in coordinate system $(x_k, x_{k+1})$ was plotted from the sequence (Fig.1).

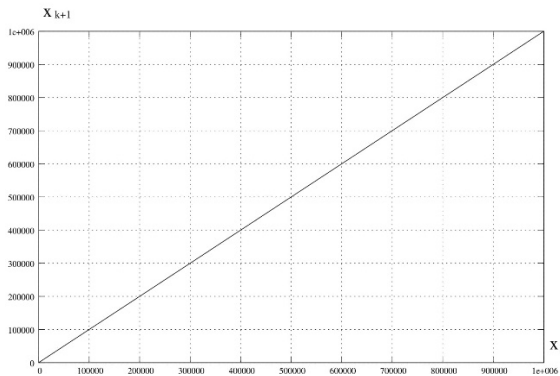

Fig.1.

The values of the *k-th* prime number were designated on the horizontal axis, whereas the values of *(k+1)-th* prime number on the vertical axis. The magnification of a fragment of the graph in Fig.1 shows that it is not a straight line inclined at the angle of 45°, but an irregular broken line positioned directly above the straight line (Fig.2). Thus, this is a certain regularity.

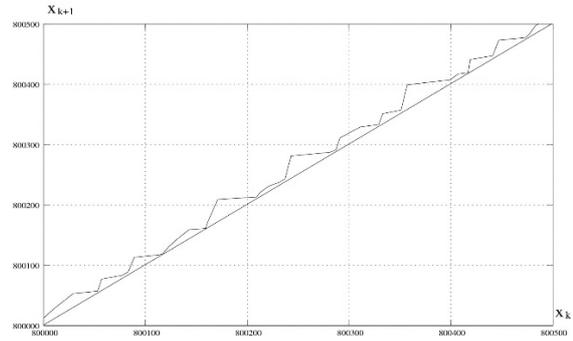

Fig. 2.

In the next step a graph which is a derivative of the broken line from Fig.2 was determined with the use of equation:

$$D_k = \frac{f(k+1)-f(k)}{f(k)-f(k-1)} \qquad (1)$$

where $f(k)$ is a function based on natural numbers, which, in this case, are prime numbers. If – for example- they were square numbers, $f(k) = k^2$ and then $D_k = \frac{2k+1}{2k-1}$, whereas, if they were even numbers, $f(k) = 2k$ and then $D_k=1$. In the case of prime numbers function $f(k)=x_k$. Therefore, an explicit mathematical form of $D_k$ is not known.

Nevertheless, a graph of $D_k$ may be designated using the above mentioned sequence of prime numbers (Fig.3).

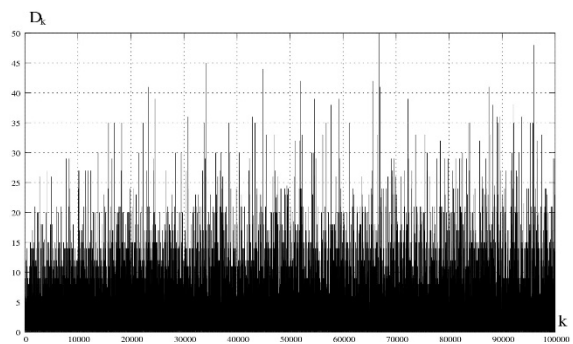

Fig. 3.

Further analysis confirmed the supposition that such graph is chaotic. A map in coordinate system ($D_k$, $D_{k+1}$) was used (Fig. 4).

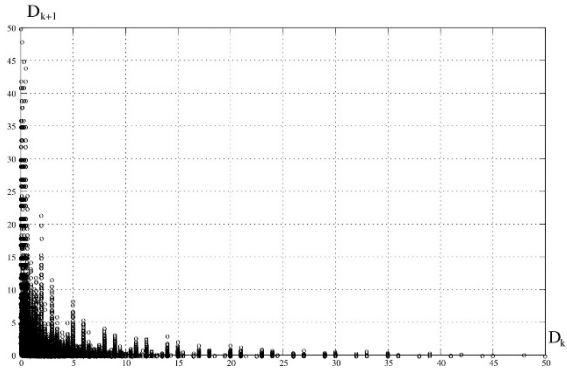

Fig. 4.

The magnification of its fragment indicates that we are dealing with a distribution that is quite regular (Fig.5).

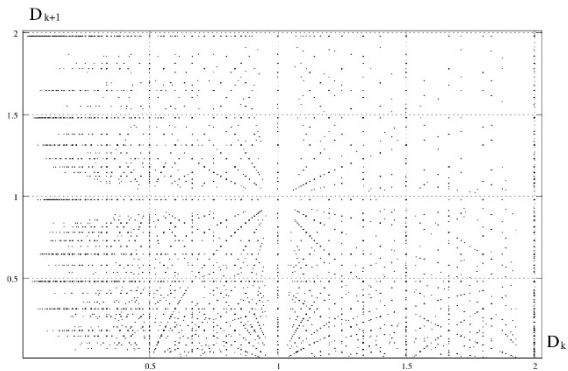

Fig. 5.

A certain regularity of the distribution of prime numbers was also discovered by Stanisław Ulam [1,2] and Martin Gardner [3]. However, it seems that the graph in Fig.5. is characterised by more regularity and much less noise than Ulam spiral. Despite this, the regularity is illusory, as it is neither a periodic distribution, nor a quasiperiodic one. It is a chaotic distribution, which is proved by the spectral analysis conducted further below. Accordingly, the amplitude spectrum was determined from the following equation:

$$Y[k] = \frac{1}{N}\sum_{n=0}^{N-1} D_n e^{-\frac{2\pi i k n}{N}} \qquad (2)$$

where: $D_n$ – $n$-th derivative sample, $i$ – imaginary unit, $N$ – number of samples, $Y[k]$ – $k$-th harmonic. As a result, the graph presented in Fig.6. was created.

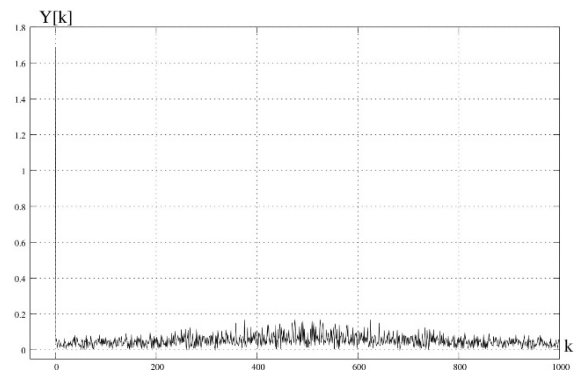

Fig. 6.

The graph shows a relatively big value of the constant component and low-amplitude noise resulting from chaos, containing a big quantity of harmonics. Theoretically, it contains an infinite number of harmonics.

Summing up, it should be stated that despite the big regularity in Fig.5., the distribution of prime numbers is chaotic. This is confirmed by the amplitude spectrum in Fig.6. In consideration of the unpredictability of chaotic phenomena, it should be claimed that it is not possible to determine a priori the position of the next prime number on the numerical axis .

Another analysis conducted in the paper concerned number π, and, specifically, the distribution of its digits after the coma. To obtain this, 50000 digits were used, their distribution on the numerical axis was presented in Fig.7.

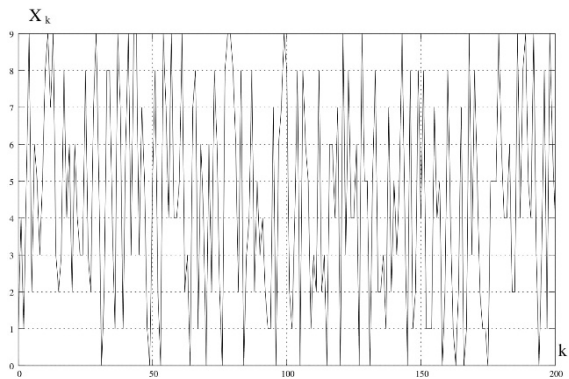

Fig. 7.

The distribution may suggest chaos. Similarly to the analysis discussed previously in the paper, a graph consisting of the digits was formed in coordinate system ($X_k$, $X_{k+1}$) (Fig.8).

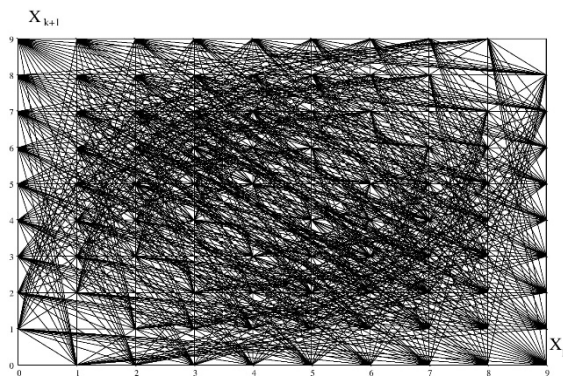

Fig. 8.

The values of the *k-th* prime number were designated on the horizontal axis, whereas the values of *(k+1)-th* prime number on the vertical axis. As a result, a fairly regular map was obtained, see Fig.8.

Just as in the previously discussed case, the regularity is illusory, because the spectral analysis of the sequence from Fig.7 indicated that we are dealing with chaos (Fig.9).

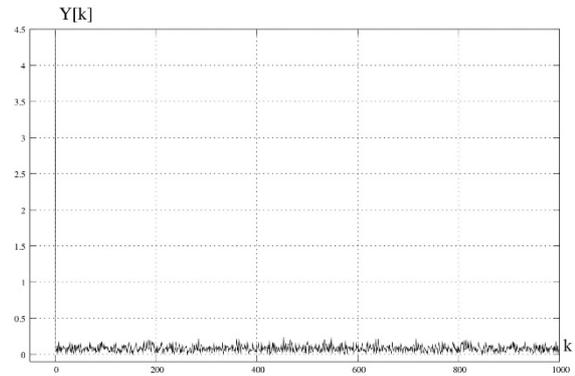

Fig. 9.

To conclude, it was proved that both the distribution of prime numbers and the distribution digits of number π have a chaotic character.